\documentclass[11pt]{article}

\usepackage{amssymb}
\usepackage{epsfig}

\newtheorem{Lemma}{Lemma}
\newtheorem{Proposition}[Lemma]{Proposition}

\newtheorem{Theorem}[Lemma]{Theorem}

\newtheorem{OP}{Open Problem}

\newcommand{\Reals}{{\mathbb{R}}}
\newcommand{\Cube}{{\mathbb{C}}}
\newcommand{\MM}{\mathcal{M}}

\newcommand{\cd}{\ \stackrel{d}{\rightarrow} \ }

\newcommand{\PP}{\mbox{${\cal P}$}}
\newcommand{\EE}{\mbox{${\cal E}$}}
\newcommand{\VV}{\mbox{${\cal V}$}}
\newcommand{\GG}{\mbox{${\cal G}$}}

\newcommand{\FF}{\mbox{${\cal F}$}}

\newcommand{\eps}{\varepsilon}

\newcommand{\var}{{\rm var}}

\newcommand{\sfrac}[2]{{\textstyle\frac{#1}{#2}}}

\newcommand{\len}{\, \mathrm{len}}

\newcommand{\pigreedy}{\pi_{\mbox{\tiny greedy}}}

\def\PI{\mathbf{\Pi}}

\newcommand{\E}{E}

\def\f{\frac}

\def\suml{\sum\limits}

\def\proof{\par\noindent{\bf Proof.\ }}
\def\proofof#1{\par\noindent{\bf Proof of #1.\ }}
\def\eop{\vskip 3mm }

\def\eqref#1{(\ref{eq.#1})}

\def\putfigure#1#2{
	\begin{figure}[ht]
	\centering
	\includegraphics{#1.eps}
	\caption{#2}
	\label{fig.#1}
	\end{figure}
}

\begin{document}

\title{Percolating paths through random points}
\author{David J. Aldous\thanks{Research supported by NSF Grant DMS-0203062}\\
Department of Statistics\\ 
367 Evans Hall \#\  3860\\
U.C. Berkeley CA 94720\\  aldous@stat.berkeley.edu
\and
Maxim Krikun\\
Department of Statistics\\ 
367 Evans Hall \#\  3860\\
U.C. Berkeley CA 94720\\  krikun@stat.berkeley.edu}

\maketitle

\begin{abstract}
We prove consistency of four different approaches to formalizing the
idea of minimum average edge-length in a path linking some
infinite subset of points of a Poisson process. 
The approaches are 
(i) shortest path from origin through some $m$ distinct 
points; 
(ii) shortest average edge-length in paths across the diagonal of a large cube;
(iii) shortest path through some specified proportion $\delta$ of points
in a large cube;
(iv) translation-invariant measures on paths in $\Reals^d$ which contain
a proportion $\delta$ of the Poisson points.
We develop basic properties of a normalized average length function
$c(\delta)$ and pose challenging open problems.
\end{abstract}

\vspace{0.1in}


\vspace{0.1in}


\newpage
\section{Introduction}
\label{sec-INT}
Fix dimension $d \geq 2$ 
and let $(\xi_i)$ be a Poisson point process of rate $1$ per unit volume in $\Reals^d$.
``Volume" means $d$-dimensional Lebesgue measure.

To start with an analogy,
one can informally describe the critical value for
continuum percolation as the smallest $c$ such that there exists some
infinite sequence 
$\xi_{j_1},\xi_{j_2},\xi_{j_3}, \ldots$
of distinct points such that
$\max_{i \geq 1} |\xi_{j_{i+1}} - \xi_{j_i}| \leq 2c$.
Here $| \ \cdot \ |$ denotes Euclidean distance.
What we study in this paper can analogously be described informally
as the smallest $c$ such that there exists some
path through an
infinite sequence 
$\xi_{j_1},\xi_{j_2},\xi_{j_3}, \ldots$
of distinct points 
whose {\em average edge-length} 
$\lim_{n \to \infty} n^{-1} 
\sum_{i = 1}^n |\xi_{j_{i+1}} - \xi_{j_i}| \leq c$.
One could formalize this directly by e.g. replacing {\em lim} by {\em lim sup}, 
requiring the property to hold almost surely and then taking the {\em inf} 
of such $c$.  
But such a definition seems neither elegant nor convenient.
Our purpose in this paper is to study four indirect approaches to this formalization question
and show that they lead to the same constant, which we call $c(0+)$.
Along the way we introduce a function $c(\delta)$ which plays a role analogous
to the percolation function.
Our results are collected in Theorem \ref{T1}.

Conceptually, this topic seems intermediate between
{\em first passage percolation} 
and the random 
{\em traveling salesman problem} (TSP).
Regarding the former, our $c(0+)$ differs 
(informally speaking)
from a continuum analog of the time constant
in first passage percolation
\cite{howard-fpp,kesten-FPP} because we use ``distance along a path" in place of 
``end-to-end Euclidean distance" 
(nonetheless we use 
``continuum first passage percolation" as a descriptor of one of our approaches below).
Regarding the latter,
take $n$ random points in a $d$-dimensional cube of volume $n$.
Let $L_n(1)$ be the length of the shortest cycle through all $n$ points,
i.e. the length of the solution of the TSP.  
Almost 50 years ago, 
Beardwood-Halton-Hammersley \cite{BHH59}
proved there exists a constant $c(1)$
such that $EL_n(1) \sim c(1)n$ as $n \to \infty$.
Subsequent work on related problems is described in the monographs
by Steele \cite{steele97} and Yukich \cite{yukich-book}.
Another of our approaches modifies 
the TSP
by considering cycles through some specified proportion $\delta$
of the Poisson points.
Additional motivation for the current paper comes from work on such problems in the ``mean-field" setting, described in section \ref{sec-mf}.

Challenging problems for future research 
are listed in section \ref{sec-discuss}. 
One of our techniques -- subadditive analysis of optimal
cost/reward ratios -- seems potentially applicable in
other contexts, as briefly discussed in section \ref{sec-newsub}.

\section{The equivalence theorem}
\label{sec-Equiv}
\paragraph{Approach 1} {\em  Continuum first passage percolation from the origin.} 
For each $m \geq 1$ define a random variable $T_m$ as the minimum, over
choices 
$\{\xi_{j_1},\ldots,\xi_{j_m}\}$ 
of $m$ distinct points of the Poisson process on $\Reals^d$, of
\[ |\xi_{j_1}|
+ \sum_{i=2}^m | \xi_{j_i} - \xi_{j_{i-1}}| . \]

\noindent
{\em Comment.}
One certainly expects that $T_m/m$ should converge to a constant,
but we see no easy argument.  
In particular we don't see how to apply subadditivity arguments
directly to study $T_m$.

\paragraph{Approach 2} {\em  Continuum first passage percolation across a diagonal.} 
For $s > 0$ define a random variable $W_s$ as the minimum, over 
all $m \geq 1$ and all choices 
$\{\xi_{j_1},\ldots,\xi_{j_{m-1}}\} \subset [0,s]^d$ 
of $m-1$ distinct points of the Poisson process, of
\begin{equation}
 m^{-1} \sum_{i=1}^{m} | \xi_{j_i} - \xi_{j_{i-1}}|  \label{mratio}
\end{equation}
where $\xi_{j_0} = (0,\ldots,0)$
and $\xi_{j_{m}} = (s,\ldots,s)$.

\noindent
{\em Comment.}
Here we can attempt subadditivity analysis,
based on splitting the cube of side $2s$ into $2^d$ subcubes
of side $s$, though because of the ``ratio" form of (\ref{mratio})
we are not in the usual format for the subadditive ergodic theorem.

\paragraph{Approach 3} {\em  TSP on sparse subsets of the cube.}
Let $\Cube_n = [0,n^{1/d}]^d$ be the cube of volume $n$ in $\Reals^d$.
Put $n$ random (independent, uniformly distributed) points
$(\zeta_j)$
into $\Cube_n$.
Fix $0<\delta \leq 1$.
Let $L_n(\delta)$ be the minimum, over all choices of cycles 
$(\zeta_{j_1}, \zeta_{j_2}, \ldots,\zeta_{j_m}, \zeta_{j_{m+1}} = \zeta_{j_1})$
through $m = \lceil\delta n \rceil$ disjoint choices from the random points, of the cycle length
$\sum_{i=1}^{m} |\zeta_{j_{i+1}} - \zeta_{j_i}|$.

\noindent
{\em Comment.}
Here subadditivity can be applied in familiar ways.
Note that, in contrast to continuum percolation where definitions
are in terms of the process on infinite space $\Reals^d$,
the approaches above all envisage taking limits over finite regions.
For the record we give a
final approach which does work directly on infinite space, 
as ``the $n = \infty$ analog of Approach 3",
though we admit it does not seem very useful.

\paragraph{Approach 4} {\em  Translation invariant distributions on infinite paths through Poisson points.}
Consider a locally finite set $(x_i)$ of points in $\Reals^d$, together 
with a set $\EE$ of edges whose endpoints are in $(x_i)$,
where the edges form a collection
of doubly-infinite paths,
each point of $(x_i)$ appearing either once or never in the paths.
Write $\mathbf{S}$ for the space of such points-and-paths configurations.
The Euclidean translation group acts naturally on $\mathbf{S}$, so one can define a probability distribution
$\mu$ on $\mathbf{S}$ to be {\em invariant} if it is invariant under the action of the
Euclidean translation group.
Let $\MM$
be the set of invariant distributions on $\mathbf{S}$
under which the distribution of the points $(\xi_i)$ is the Poisson point process of rate $1$.
Informally, a $\mu \in \MM$ is just a way of collecting some subsets of the Poisson points into
paths using a rule which doesn't depend on the location of the origin.
For $\mu \in \MM$ there is a constant 
$\delta(\mu) \in [0,1]$
specified informally as
``the proportion of points which are in some infinite path"
and defined formally via the formula: for every cube $\Cube \subset \Reals^d$,
\[ E_\mu (\mbox{number of points $\xi_i \in \Cube$ 
which are in some infinite path}) = 
\delta(\mu) \ \mbox{volume}(\Cube) .\]
Similarly there is a constant $\ell(\mu)$ interpreted as
``mean edge-length over all edges in the paths of $\mu$"
and formally via the formula: for every cube $\Cube \subset \Reals^d$,
\[ E_\mu(
\mbox{length of } \EE \cap \Cube)
=
\delta(\mu) \ell(\mu)\ \mbox{volume}(\Cube) \]
(here we regard $\EE$ as a subset of $\Reals^d$).
Finally define
\[
\bar{c}(\delta) := \inf \{\ell(\mu): \ \mu \in \MM, \ \delta(\mu) = \delta\}
. \]

\begin{Theorem}~
\label{T1}
\begin{itemize}
\item[(a)] 
  For $0<\delta \leq 1$ there exists a constant $c(\delta)$ such that
  $\frac{L_n(\delta)}{\delta n} \to c(\delta)$ in $L^2$ as $n \to \infty$. 
\item[(b)] 
  The function $c(\delta)$ is non-decreasing and continuous on $(0,1]$,  
  the function $\delta c(\delta)$ is convex,
  and the limit $c(0+) := \lim_{\delta \downarrow 0} c(\delta)$
  is strictly positive.
\item[(c)] $W_s \to c(0+)$ a.s. as $s \to \infty$.
\item[(d)] $\bar{c}(\delta) = c(\delta), \ 0 < \delta \leq 1$.
\item[(e)] $m^{-1}T_m \to c(0+)$ in probability.
\end{itemize}
\end{Theorem}
From (a) we see $c(1)$ is the constant in the
Beardwood-Halton-Hammersley \cite{BHH59}
theorem.
Monte Carlo simulations \cite{JRS04} give 
(for $d = 2$)
$c(1) \approx 0.7119$
but no close rigorous bounds are known.

\subsection{Discussion}
\label{sec-discuss}
The function $c(\delta), \ 0 < \delta \leq 1$
seems worthy of study as a analog of the classical
percolation function
from lattice percolation theory \cite{gri99}:
\[ f(p) := P(
\mbox{origin is in some infinite component}) \]
in bond percolation with edge-probability $p$.
In particular one can ask whether there exists a
{\em scaling exponent} $0< \alpha < \infty$,
that is whether
\begin{equation}
c(\delta) - c(0+) \asymp \delta^\alpha 
\mbox{ as } \delta \downarrow 0 . \label{def-alpha}
\end{equation}
For the record we state the (probably very hard)
\begin{OP}
Prove (\ref{def-alpha}) holds for some $\alpha$.
Or give bounds on the possible values of $\alpha$.
\end{OP}
The next question may be easier.
The facts that $c(\delta)$ is nondecreasing and that
$\delta c(\delta)$ is convex imply: either\\
(i) $c(\delta)$ is strictly increasing on $0<\delta<1$; or\\
(ii) $c(\delta)$ is constant on $0<\delta<\delta_0$,
for some $0<\delta_0 \leq 1$.
\\
But the latter seems implausible.
\begin{OP}
Prove that $c(\delta)$ is strictly increasing on $0<\delta<1$.
\end{OP}
Getting reasonable bounds on the numerical value of $c(0+)$ 
seems difficult.  
Standard methods
(comparison with branching random walk: section \ref{sec-cLB})
give an explicit lower bound (\ref{c0+LB}), which in $d=2$ is 
$1/(2 \pi e) \approx 0.0585$.
But we don't see any simple way to get an interesting upper bound.
Even Monte Carlo methods seem difficult to code convincingly;
for the record we write
\begin{OP}
In $d=2$ study the numerical value $c(0+)$ and the presumed scaling
exponent $\alpha$ via Monte Carlo methods.
\end{OP}
The variation in which 
(in $d=2$)
we restrict to ``upward oriented" paths,
that is edges $(x_i,y_i) \to (x_{i+1},y_{i+1})$
are required to have $y_{i+1} > y_i$,
is easier to study via simulation; our small-scale
simulations suggest the analog of $c(0+)$ in this variation is
$\approx 0.62$, which would be an {\em a priori} upper bound
for our original $c(0+)$.

Another question concerns variances.
Take $d = 2$ here.
For the basic traveling salesman problem, 
that is for $L_n(1)$,
it is known that 
$\var (L_n(1)) $ is order $n$; precisely, 
\[ 
\liminf_n n^{-1} \var (L_n(1)) > 0, \quad 
\limsup_n n^{-1} \var (L_n(1)) < \infty . \]
The upper bound is explicit in Steele \cite{ste81}
and the lower bound follows from a corresponding large deviation
lower bound in 
Rhee \cite{rhee-91}.
On the other hand, for first-passage percolation it 
has long been conjectured \cite{durrett-kesten} that variance grows as the
$\frac{2}{3}$ power of expectation,
though little has been proved rigorously \cite{BKS03}.
\begin{OP}
Prove 
$\var (T_m) \asymp m^{2/3}$.
Or just prove
$\var (T_m) = o(m)$.
\end{OP}
We will not even venture a conjecture for the asymptotic behavior
of $\var (L_n(\delta))$.

As a small rewriting of the definition of $L_n(\delta)$, 
let $L(n,m)$ be the minimum, over all choices of cycles 
$(\zeta_{j_1}, \ldots,\zeta_{j_m}, \zeta_{j_{m+1}} = \zeta_{j_1})$
through some chosen $m $ of the random points, of the cycle length
$\sum_{i=1}^{m} |\zeta_{j_{i+1}} - \zeta_{j_i}|$.
Theorem \ref{T1} implies that, if $m_n/n \to 0$ sufficiently slowly,
then 
$L(n,m_n)/m_n \to c(0+)$.
However it seems plausible this also holds for smaller values of $m_n$.
\begin{OP}
Prove that 
$L(n,m_n)/m_n \to c(0+)$
in probability whenever $m_n/n \to 0$ and
$\frac{m_n}{\log^a n} \to \infty \ \forall a < \infty$.
\end{OP}
In other words, for fixed $c < c(0+)$ consider Poisson points in a cube of volume $n$;
is it true that any cycle with average edge-length $\leq c$ can have at most 
poly-log$(n)$ edges?
This would be an analog of the fact that subcritical percolation cluster size distribution has a
geometrically-decreasing tail 
\cite{gri99}.

\subsection{The mean-field model}
\label{sec-mf}
Instead of the Euclidean model in this paper,
one can consider a ``mean-field" model 
on $n$ points for which the ${n \choose 2}$ inter-point links
are assumed to have independent random lengths with Exponential (mean $n$)
distribution.
Within this model
one can define
the function $\tilde{c}(\delta)$ analogous to $c(\delta)$.
M{\'e}zard-Parisi \cite{MP86} used the non-rigorous 
{\em replica method} of statistical physics to argue
$n^{-1} EL_n(1) \to \tilde{c}(1)$,
where $\tilde{c}(1) \approx 2.04$ is derived from numerical solution of a certain
fixed point equation.
Using probabilistic reformulations of these statistical physics ideas,
Aldous \cite{me109} gave a (still non-rigorous) analysis  of the
whole function $\tilde{c}(\delta)$, 
exhibited in Figure 1 of \cite{me109},
which suggests
$\tilde{c}(\delta) - \tilde{c}(0+) \asymp \delta^{1/3}$
as $\delta \downarrow 0$.
So one might 
conjecture that the scaling exponent $\frac{1}{3}$
also holds in the Euclidean case.  

Note that the subadditivity arguments
we use in the Euclidean case to prove Theorem \ref{T1} rest upon the ``boundary effects are negligible"
property
of $\Reals^d$.
In the mean-field model, the limit analog of the Poisson process
is a certain infinite random tree, for which boundary effects are
not negligible and subadditivity arguments cannot be used.
Indeed, Approach 4 was developed in the mean-field setting as a
substitute for subadditivity.

\section{Proofs}
\label{sec-proofs}
We start with Approach 3 and prove part (a) of Theorem~\ref{T1} 
in \ref{sec-sparse};
then we prove part (b) of Theorem~\ref{T1} in \ref{sec-cdelta} 
and give a lower bound on $c(0+)$ in \ref{sec-cLB}.
Sections \ref{sec-Ws}, \ref{sec-tid}, \ref{sec-Tm} contain
the proofs of the remaining parts (c), (d), (e) of Theorem~\ref{T1}.
Part (e) seems hardest, for reasons explained at the start of section 
\ref{sec-Tm}.

\subsection{TSP on sparse subsets} \label{sec-sparse}
We rely on subadditivity arguments as in \cite{steele97,yukich-book}.
These monographs develop general results for sub- or superadditive
Euclidean functions satisfying regularity properties.
Unfortunately functionals like $L_n(\delta)$ lack the
{\em monotonicity property} (\cite{steele97} equation (3.5))
and it is not clear whether the {\em smoothness property}
(\cite{yukich-book} section 3.3) is both valid and exploitable.
We will use the inequality in (b) below as a substitute
for monotonicity.

Write $A_1,A_2,\ldots$ for constants depending only on
dimension $d \geq 2$.
We start with a purely deterministic lemma
(note that by scaling the case of general $s$ is equivalent to
the case $s=1$).
\begin{Lemma}\label{L1}
Let $\{x_1,\ldots,x_n\}$ be arbitrary points in the cube $[0,s]^d$
and let $L(m)$ be the length of the shortest cycle through some 
$m\le n$ of these points.
\begin{itemize}
\item[(a)] {\bf [Uniform boundedness]} 
    \begin{equation}\label{eq.unibo}
    L(n) \leq A_1 s n^{(d-1)/d} .
    \end{equation}
\item[(b)] For $1 \leq m_1 < m_2 \leq n$,
\[ \frac{L(m_1)}{m_1} 
   \leq \frac{L(m_2)}{m_2} + \frac{sd^{1/2}}{m_1}   .
\]
\item[(c)] {\bf [Geometric subadditivity]} 
Let $k \geq 2$ and let $(\Cube^j, \ 1 \leq j \leq k^d)$
be the natural partition of $[0,s]^d$ into $k^d$ subcubes
of side $s/k$.
Suppose that, for each $j$, there exists a cycle of length
$l_j$ through some subset 
$S_j \subset \{x_i\} \cap \Cube^j$.
Then there exists a cycle through $\cup_j S_j$ of length at most
$\sum_j l_j  + A_2 s k^{d-1}$.
\end{itemize}
\end{Lemma}

\proof 
Parts (a) and (c) are standard (\cite{steele97} sections 2.2 and 2.3).
For (b), let $y_1, y_2, \ldots, y_{m_2}, y_{m_2+1} = y_1$
be a minimum-length cycle attaining $L(m_2)$.
Then there exists $k$ such that (interpreting $k+i$ modulo $m_2$)
\[ \frac{1}{m_1} \sum_{i=1}^{m_1} 
   |y_{k+i+1} - y_{k+i}| \leq \frac{L(m_2)}{m_2} 
\]
because the right side equals the {\em average} of the left side
as $k$ varies.
To make a cycle on 
$\{y_{k},\ldots,y_{k+m_1 - 1}\}$
replace edge $(y_{k+m_1-1},y_{k+m_1})$
by edge $(y_{k+m_1-1},y_k)$, whose length is at most the diameter
$s d^{1/2}$ of the cube $[0,s]^d$.
\eop

{\em Remark.}
Lemma \ref{L1}(a) implies that the worst-case cycle length is the same order of magnitude
as the average-case lengths we will be studying.
This has the pleasant consequence that events of probability tending to zero will
make asymptotically negligible contributions to expectation of length, and so
can be ignored: we use this
{\em uniform boundedness property} several times later.

We start analysis of the probability model by making definitions to which
subadditivity arguments can easily be applied.
Recall $(\xi_i)$ denotes a Poisson point process of rate $1$ per unit volume in $\Reals^d$.
Let $N(s)$ be the number of points of $(\xi_i)$ in
$[0,s]^d$.
Define $L(s,\delta)$ as the length of the minimum-length cycle through at least
some $\lceil \delta N(s) \rceil$ of the points of
$(\xi_i) \cap [0,s]^d$.
By the triangle inequality, replacing 
``at least 
$\lceil \delta N(s) \rceil$" 
by ``exactly 
$\lceil \delta N(s) \rceil$ points"   
changes nothing.
We will start by using the next lemma as a definition of $c(\delta)$, 
and later show this agrees with the limit in Theorem \ref{T1}(a).

\begin{Lemma}
\label{L2}
For fixed $0<\delta \leq 1$
there exists a constant $c(\delta)$ such that
\[ \frac{L(s,\delta)}{\delta s^d} \to c(\delta) \mbox{ in $L^2$ as } s \to \infty. 
\]
\end{Lemma}
\proof
First note 
$s \to EL(s,\delta)$ is continuous, 
from the representation
\[ \E L(s,\delta) 
   = \sum_{n=0}^{\infty} \f{s^{nd}}{n!} e^{-s^d} \, s a_n(\delta),
\]
where $a_n(\delta)$ is the expected length of the shortest cycle through some
$[\delta n]$ of $n$ uniform random points in the unit cube $[0,1]^d$.
Next, given $s>0$, we can write any $x>0$ as $x=ks+t$ with integer $k$ and $t\in [0,s)$.
Geometric subadditivity (Lemma \ref{L1}(c)) then implies
\[ \E L(x,\delta) 
   \leq k^d \E L(s+t/k,\delta) + A_2 (s+t/k)k^{d-1}.
\]
Taking $x\to\infty$ while keeping $s$ fixed  
\[
   \limsup_{x\to\infty} x^{-d} \E L(x,\delta)
   \leq 
   \inf_{\eps>0} \, \sup_{\tau \in [0,\eps]} s^{-d} \E L(s+\tau,\delta).
\]
So by the continuity property
\begin{equation}
 \lim_{s\to\infty} s^{-d} \E L(s,\delta)
     = \inf_{s>0} s^{-d} \E L(s,\delta) 
     = \delta c(\delta) 
, \mbox{ say}
\label{qwe}
\end{equation}
with $0 \leq c(\delta) < \infty$.


Once again by geometric subadditivity,
temporarily abbreviating $L(s,\delta)$ to $L_s$,
\[ \E L_{ks}^2 \le k^d \E L_s^2  +  k^d (k^d-1) (\E L_s)^2 
                   + A_2 s k^{2d-1} \E L_s + (A_2 s k^{d-1})^2.
\]
By the uniform boundedness property \eqref{unibo}
\[ \E L_{s}^2 \le A_1^2 s^2 \E \PP_{(s^d)}^{\f{2(d-1)}{d}}
\le A_1^2s^2 \cdot A_3s^{2(d-1)}
   = A_1^2 A_3 s^6,
\]
where $\PP_{(s^d)}$ has Poisson$(s^d)$ distribution
and where the second inequality holds for some $A_3$ for all $s \geq 1$.
Combining the two displayed inequalities and (\ref{qwe}) gives
\[ \limsup_x x^{-2d} \E L_x^2 = (\delta c(\delta))^2, \]
and finally
\[ s^{-d} L(s,\delta) \to \delta c(\delta) \mbox{ in } L^2 . \]
\eop

{\em Remark.}
More sophisticated modern proofs
(\cite{steele97} sec. 2.4; \cite{yukich-book} sec. 4.1)
of the TSP case ($\delta = 1$)
use concentration inequalities
to obtain almost sure convergence; we have not investigated
concentration inequalities or a.s. convergence for $L(s,\delta)$.

%

Now we proceed with the proof of Theorem \ref{T1}(a)
with $c(\delta)$ defined by Lemma \ref{L2}.

\proofof{Theorem~\ref{T1}(a)}
Recall the definition of $L_n(\delta)$.
Let $\Cube_n = [0,n^{1/d}]^d$ be the cube of volume $n$ in $\Reals^d$.
Put $n$ random (independent, uniformly distributed) points
$(\zeta_j)$
into $\Cube_n$.
Fix $0<\delta \leq 1$.
Let $L_n(\delta)$ be the minimum, over all choices of cycles 
$(\zeta_{j_1}, \zeta_{j_2}, \ldots,\zeta_{j_m}, \zeta_{j_{m+1}} = \zeta_{j_1})$
through any chosen $m = \lceil\delta n \rceil$ of the random points, of the cycle length
$\sum_{i=1}^{m} |\zeta_{j_{i+1}} - \zeta_{j_i}|$.
Again by the triangle inequality, this is the same as
saying 
``any $m \geq \lceil\delta n \rceil$ of the random points".  

Fix small $\eps > 0$ and consider, in $\Cube_n$, a Poisson
process of rate $1 - \eps$ per unit volume.
By standard properties of the Poisson process, for each $n$ we can
couple this to the process of $n$ i.i.d. uniform points
in $\Cube_n$ in such a way that, with probability
$\to 1$ as $n \to \infty$, 
each point of the Poisson process is a point of the uniform process.
Call this the {\em inclusion coupling}.
There is a similar inclusion coupling between the uniform
process on $\Cube_n$ and the Poisson process of rate $1+ \eps$.

Now write $L_\lambda (s,\delta)$ 
to mean the quantity $L(s,\delta)$ applied
to a Poisson process of rate $\lambda$;
and write $N_\lambda(s)$ for the number of points of that
Poisson process in $[0,s]^d$.
When the inclusion couplings hold 
and when 
\[ \frac{\delta}{1+2\eps} N_{1+\eps}(n^{1/d})
\leq \delta n
\leq \frac{\delta}{1-2\eps} N_{1-\eps}(n^{1/d})
\]
then we have
\begin{equation}
 L_{1+\eps}(n^{1/d},\sfrac{\delta}{1+2\eps})
\leq L_n(\delta)
\leq L_{1-\eps}(n^{1/d},\sfrac{\delta}{1-2\eps})
 . \label{LLL}
\end{equation}
This holds because, in informal language, for each $\leq$ we
have more points to choose from, and a weaker constraint on 
minimum number of points in the cycle.
Now by scaling
\[ L_{1-\eps}(n^{1/d},\sfrac{\delta}{1-2\eps})
   = (1-\eps)^{-1/d} 
     L((1-\eps)^{-1/d}n^{1/d}, \sfrac{\delta}{1-2\eps}), 
\]
and similarly 
\[ L_{1+\eps}(n^{1/d},\sfrac{\delta}{1+2\eps})
   = (1+\eps)^{-1/d} 
     L((1+\eps)^{-1/d}n^{1/d}, \sfrac{\delta}{1+2\eps}).
\]
Taking limits in (\ref{LLL}), using Lemma \ref{L2} and
the uniform boundedness property, we get
\[ \limsup_n n^{-1} E L_n(\delta) 
   \leq (1-\eps)^{-1-1/d} 
   \sfrac{\delta}{1-2\eps} c(\sfrac{\delta}{1-2\eps}) 
\]
and similarly
\[ (1+\eps)^{-1-1/d} 
\sfrac{\delta}{1+2\eps}
c(
\sfrac{\delta}{1+2\eps}
) 
\leq
 \liminf_n n^{-1} E L_n(\delta) 
. \]
Letting $\eps \downarrow 0$ and using continuity of $c(\delta)$ 
(which we prove independently in the next subsection)
we see that for $0<\delta <1$
\[\frac{L_n(\delta)}{\delta n} \to c(\delta) \mbox{ in }L^1.  \] 
The case $\delta = 1$ is similar, but of course is already part of the usual proof
\cite{steele97} of the Beardwood-Halton-Hammersley theorem,  
so we omit it.

\subsection{Properties of $c(\delta)$} \label{sec-cdelta}
Next we prove part (b) of Theorem \ref{T1}.

\begin{Proposition}$c(\delta)$ is non-decreasing
on $0<\delta \leq 1$.
\end{Proposition}
\proof 
Fix $0<\delta_1 < \delta_2 \leq 1$.
By Lemma \ref{L1}(b)
\[ \frac{L(s,\delta_1)}{\lceil \delta_1 N(s) \rceil}
   \leq \frac{L(s,\delta_2)}{\lceil \delta_2 N(s) \rceil}
    + \frac{sd^{1/2}}{\lceil \delta_1 N(s) \rceil} \]
and so
\[
   \frac{L(s,\delta_1)}{\delta_1 s^d}
   \leq R_s \frac{L(s,\delta_2)}{\delta_2 s^d}
   + \frac{sd^{1/2}}{\delta_1 s^d}
\]
where
\[ R_s = \frac
    {\lceil \delta_1 N(s) \rceil /\delta_1}
    {\lceil \delta_2 N(s) \rceil /\delta_2}.
\]
Since $R_s$ is uniformly bounded and $R_s \to1$ as $s\to\infty$,
using Lemma~\ref{L2} we deduce $c(\delta_1) \leq c(\delta_2)$.
\begin{Proposition} $\delta c(\delta)$ is convex 
on $0<\delta \leq 1$.
\end{Proposition}
\proof 
Fix $0<\delta_1 < \delta_2 \leq 1$ and $0<\lambda<1$.
Let $(\Cube^j, \ 1 \leq j \leq k^d)$
be the partition of $[0,s]^d$ into $k^d$ equal subcubes.
Let $N_j$ be the number of points of $(\xi_i)$ in the $j$-th subcube and
let $\ell_j(s)$ (resp. $\tilde\ell_j(s)$) be the length of the shortest
cycle through some $\lceil\delta_1 N_j\rceil$ (resp. $\lceil\delta_2 N_j\rceil$)
points in $\Cube^j$.

Take any $\delta  < \lambda \delta_1 + (1-\lambda)\delta_2$.
The event 
\[
  \sum_{j=1}^{\lceil \lambda k^d \rceil} \lceil \delta_1 N_j \rceil
  + \sum_{j=\lceil \lambda k^d \rceil + 1}^{k^d} \lceil \delta_2 N_j \rceil
\   \geq \  \lceil \delta N(s) \rceil
\]
has probability $\to 1$ as $s \to \infty$,
and on this event we have by Lemma \ref{L1}(c)
\[ L(s,\delta) 
   \leq \sum_{j=1}^{\lceil \lambda k^d \rceil} \ell_j(s)
        + \sum_{j=\lceil \lambda k^d \rceil + 1}^{k^d} \tilde \ell_j(s)
        + A_2 s k^{d-1}. 
\]
Taking expectations, letting $s \to \infty$ and using Lemma \ref{L2} we obtain
\[ \delta c(\delta) \leq 
   \frac{ \lceil \lambda k^d \rceil}{k^d} \ \delta_1 c(\delta_1)
    + \frac{ k^d -  \lceil \lambda k^d \rceil}{k^d} \ \delta_2 c(\delta_2),
\]
and letting $k \to \infty$
\begin{equation}
 \delta c(\delta) \leq \lambda \delta_1 c(\delta_1) 
                         + (1-\lambda) \delta_2 c(\delta_2).
\label{dcd}
\end{equation}
If $\delta = \lambda^\prime \delta_1 + (1-\lambda^\prime)\delta_2$
then (\ref{dcd}) holds for all
$\lambda < \lambda^\prime$
and hence for
$\lambda = \lambda^\prime$, 
proving convexity.

\begin{Proposition} $c(\delta)$ is continuous on $(0,1]$. \end{Proposition}
\proof 
Convexity of $\delta c(\delta)$ implies continuity of $\delta c(\delta)$,
and hence continuity of $c(\delta)$, on the open interval $0<\delta<1$.
Continuity at $\delta = 1$ requires a separate argument.

Take $\delta < 1$. The cycle attaining $L_n(\delta)$ passes through
$\lceil n \delta \rceil$ points in $[0,n^{1/d}]^d$.
By Lemma \ref{L1}(a) there exists a cycle through the
remaining $q = n - \lceil n \delta \rceil$ points
with length at most
\[ A_1 n^{1/d} q^{(d-1)/d} \leq A_1 n (1- \delta)^{(d-1)/d} . \]
By joining the two cycles we find
\[ L_n(1) \leq L_n(\delta) 
+ A_1 n (1-\delta)^{(d-1)/d}
+ 2 n^{1/d} d^{1/2} . \]
Letting $n \to \infty$ 
\[ c(1) \leq \delta c(\delta) + A_1 (1-\delta)^{(d-1)/d} \]
and this implies continuity as $\delta \uparrow 1$.



\subsection{A lower bound on $c(0+)$} \label{sec-cLB}
We start by noting a one-sided bound relating the $T_k$
in Approach 1 to $c(0+)$.
In the context of Theorem \ref{T1}(a), the cycle attaining
length $L_n(\delta)$ can be converted into a path from the origin
by replacing some edge 
$(\zeta,\zeta^\prime)$
by the edge from the origin to $\zeta^\prime$.
It follows that
\[
\lim_k P( k^{-1}T_k \leq c(\delta) + \eps) = 1 \ \forall \eps > 0
\]
and thus 
\begin{equation}
\lim_k P( k^{-1}T_k \leq c(0+) + \eps) = 1 \ \forall \eps > 0
. \label{4bound}
\end{equation}
We can now use a standard argument.
Consider 
branching random walk 
on $\Reals^+$,
starting with one individual at the origin in generation $0$,
and where each individual in each generation 
(at position $x$, say)
has children at positions
$(x + |\xi_j|, \ j \geq 1)$
where $(\xi_j)$ forms a Poisson point process of rate $1$
in $\Reals^d$.
Write $\theta_k(\cdot)$ for the mean measure
for the positions of the generation-$k$ individuals
$(Y_{k,i}, i \geq 1)$:
\[ \theta_k(\cdot) = \sum_i P(Y_{k,i} \in \cdot) . \]
This is exactly the same measure as the mean measure for
lengths of $k$-step paths from the origin
through the Poisson points:
\[ \theta_k(\cdot) = 
\sum_{(j_1,\ldots,j_k)} 
P\left(|\xi_{j_1}| + \sum_{i=2}^k |\xi_{j_i} - \xi_{j_{i-1}}| \ \in \cdot
\right)
\]
where the sum is over ordered distinct $k$-tuples.
So for $T_k$ as defined in Approach 1,
and for $c>0, \lambda>0$
\begin{eqnarray*}
P(T_k \leq ck)
&\leq&
\theta_k[0,ck] \quad  \mbox{ (Markov's inequality)}\\
&\leq&
e^{\lambda ck}
\int_0^\infty e^{-\lambda x} \theta_k(dx)
 \quad \mbox{ (large deviation inequality)}\\
&=& \left[
e^{\lambda c} \int_0^\infty e^{-\lambda x} \theta_1(dx) \right]^k
\end{eqnarray*}
by the structure of branching random walk.
Comparing with (\ref{4bound}) we see
\[ c(0+) \geq \sup\left\{c: \ 
\inf_{\lambda >0} 
e^{\lambda c} 
\int_{\Reals^d} e^{-\lambda |y|} \ dy \ < 1 \right\} . \]
Writing $v_d$ for the volume of the unit ball in $\Reals^d$,
\[ 
\int_{\Reals^d} e^{-\lambda |y|} \ dy 
= \frac{v_d \Gamma(d+1)}{\lambda^{d-1}} . \]
The $\inf_{\lambda > 0} ( \ \cdot \ )$
is now attained at
$\lambda = (d-1)/c$
and we finally find
\begin{equation}
c(0+) \geq 
e^{-1} (d-1) 
\left(v_d \Gamma(d+1) \right)^{-\sfrac{1}{d-1}} 
\label{c0+LB}
\end{equation}
and of course
$v_d = \pi^{d/2}/\Gamma(1+d/2)$.
In particular, for $d=2$ we find
$c(0+) \geq (2 \pi e)^{-1}$.

\subsection{The limit for $W_s$} \label{sec-Ws}
We now turn to Approach 2.
Here we consider the rate-$1$ Poisson process $(\xi_i)$ on $\Reals^d$.
Let $\PI_{0,s}$ be the set of paths $\pi$ across the diagonal
of $\Cube_s$; that is, of paths
\[ (0,\ldots,0) = \xi_{j_0}, \xi_{j_1}, \ldots ,
\xi_{j_{m-1}}, \xi_{j_m} = (s,\ldots,s) \]
where $\{\xi_{j_1},\ldots,\xi_{j_{m-1}}\}$
are distinct points of 
$\{\xi_i\} \cap [0,s]^d$.
Write
\[ \ell(\pi) = \sum_{i=1}^{m} | \xi_{j_i} - \xi_{j_{i-1}}|  ;
   \quad 
   m(\pi) = m 
\]
for the length and number of edges in the path $\pi$
and 
\[ w(\pi) = \ell(\pi)/m(\pi) \]
for the {\em average edge-length} of $\pi$.

By considering diagonally-adjacent unit cubes and picking
(where possible) one point from each, we see there exists a path
$\pigreedy \in \PI_{0,s}$ 
such that
\begin{equation}
  \begin{array}{c}
  m(\pigreedy) - 1 \mbox{ has Binomial $(\lfloor s \rfloor, 1-e^{-1})$ distribution;}
    \vphantom{\int\limits_I}\\
  \ell(\pigreedy) \leq A_4 \lceil s \rceil
  \end{array}
\label{specialpi}
\end{equation}
for some constant $A_4$.
Recall the definition
\[ W_s := \min_{\pi \in \PI_{0,s}} \frac{\ell(\pi)}{m(\pi)}. \]
Applying to $\pigreedy$ and using (\ref{specialpi}) we see
\begin{equation}
\limsup_s W_s \leq A_4/(1-e^{-1}) \mbox{ a.s.}
\label{Ww}
\end{equation}
Similarly one can check that
$s/m(\pigreedy)$ is uniformly integrable as $s \to \infty$ and so
\begin{equation}
(W_s, 0<s<\infty)
\mbox{ is uniformly integrable}.
\label{UI}
\end{equation}
While subadditivity is not applicable directly to $W_s$,
the proof below is an easy indirect application.

\begin{Proposition}\label{C1}
$W_s \to \beta$ a.s. and in $L^1$ as $s \to \infty$,
for some constant $0 \leq \beta < \infty$.
\end{Proposition}
\proof
For $c \geq 0$ define
\[ X^{(c)}_{0,s} 
   =  \min_{\pi \in \PI_{0,s}} ( \ell(\pi) - c m(\pi) + c). 
\]
We shall prove the Proposition for
\[ \beta:= \sup \{c: E X^{(c)}_{0,s} \geq 0 \ \forall s > 0 \}. \]
For $s<t$ write $\PI_{s,t}$ for the set of paths across the diagonal
of $[s,t]^d$ and define $X^{(c)}_{s,t}$ analogously to $X^{(c)}_{0,s}$.
Given a path
$\pi_1 \in \PI_{0,s}$
and a path
$\pi_2 \in \PI_{s,t}$,
their concatenation gives a path 
$\pi \in \PI_{0,t}$,
for which
\begin{equation}
 \ell(\pi) \leq \ell(\pi_1) + \ell(\pi_2); \quad 
 m(\pi) = m(\pi_1) + m(\pi_2) -1. \label{concat}
\end{equation}
In such a concatenation, the last edge of $\pi_1$ and the first edge of $\pi_2$
are replaced by a single edge
and the inequality for $\ell(\pi)$ arises only from the triangle inequality 
for this replacement.

Consider first a value $c$ such that
$EX^{(c)}_{0,s} < 0$ for some $s$.
For this $s$  
let $\pi$ be the random path 
in $\PI_{0,s}$
such that
$E(\ell(\pi) - c m(\pi) + c)  = EX^{(c)}_{0,s} < 0$.
So by the concatenation property (\ref{concat})
and the strong law of large numbers we can construct
random paths $\pi_k$ in $\PI_{0,ks}$
such that
\[ \limsup_{k \to \infty} \f{\ell(\pi_k)}{m(\pi_k)-1}
   \leq \f{E \ell(\pi)}{E m(\pi)-1} 
   < c 
   \mbox{ a.s. } 
\]
and it easily follows that
\[ \limsup_{s \to \infty} W_s < c 
   \mbox{ a.s. } 
\]
From choice of $c$ and definition of $\beta$ we deduce
\[ \limsup_{s \to \infty} W_s \leq \beta 
\mbox{ a.s. } \]
On the other hand consider a value of $c$ such that
$EX^{(c)}_{0,s} \geq 0$ for all $s > 0$.
The process $(X^{(c)}_{s,t})$ is subadditive by (\ref{concat}).
A routine application of the
subadditive ergodic theorem
shows that there exists a constant
$\gamma(c) \geq 0$ such that
\begin{equation}
 \lim_{s \to \infty} s^{-1} X^{(c)}_{0,s} = \gamma(c)
\mbox{ a.s. } \label{sX}
\end{equation}
Now $W_s = \ell(\pi_s)/m(\pi_s)$ for some random $\pi_s$ in $\PI_s$,
and
\[ X^{(c)}_{0,s} \leq \ell(\pi_s) - c m(\pi_s) + c , \]
implying
\begin{equation}
 W_s - c \geq \frac{X^{(c)}_{0,s} - c }{m(\pi_s)} . \label{Wsc}
\end{equation}
Combining (\ref{sX}) with (\ref{Wsc}) gives
\[ \liminf_s ( W_s - c )
   \ge  \gamma(c) \cdot \liminf_s \f{s}{m(\pi_s)} \ge 0 \mbox{ a.s. }
\]
and consequently
\[ \liminf_s W_s \geq c \mbox{ a.s. } \]
From choice of $c$ and definition of $\beta$ we deduce
\[ \liminf_{s \to \infty} W_s \geq \beta 
   \mbox{ a.s. } 
\]
completing the proof of a.s. convergence in Proposition \ref{C1}.  
Finally, note $\beta < \infty$ by (\ref{Ww}),
and note that $L^1$ convergence follows from (\ref{UI}).

We use concatenation constructions based on (\ref{concat}) several times later,
and here is a rather technical formulation of the results of such 
constructions, designed to replace the use of the expectations
$E \ell(\pi_s)$ and $E m(\pi_s)$ by using truncation.
By ``adjacent cubes" we mean disjoint cubes in which we can choose
diagonals to form a connected path (e.g. as in Figure 1 below).

\begin{Lemma}[Concatenation argument]\label{Lconcat}
Let $\pi_s$ be a random path across the diagonal of the cube $[0,s]^d$,
and let $w_0\in(0,s)$.   Then there exist paths $\pi_{sn}$ through the Poisson points 
in $n$ adjacent cubes of side $s$ each, such that
\begin{equation}\label{eq.concat}
  \limsup_n w(\pi_{sn}) 
  \le w_0 + \f{ w_0 + sd^{1/2} P( w(\pi_s)>w_0 )}
            { (sd^{1/2}/w_0 - 1) P(w(\pi_s) \leq w_0)}
      \quad\mbox{a.s.}
\end{equation}
and
\[ \lim_n \f{m(\pi_{sn})}{n}
   \ge (sd^{1/2}/w_0-1) P( w(\pi_s)\le w_0 )
   \quad\mbox{a.s.} \]
\end{Lemma}
\proof
Given a random path $\pi_s$ across the diagonal of the cube $[0,s]^d$,
consider the modified path
\[ \tilde\pi_s = \cases{
	 	 \pi_s & on the event $w(\pi_s) \le  w_0 $,\cr
		 \nearrow_s & otherwise,
 		}
\]
where $\nearrow_s$ is a shortcut path consisting of a single edge
from $(0,\ldots,0)$ to $(s,\ldots,s)$.
%
By concatenating $n$ independent copies of $\tilde\pi_s$ 
we get a path $\pi_{ns}$ such that, using (\ref{concat}),
\begin{eqnarray}
 w(\pi_{ns}) 
        &\le& \f{ \suml_{i=1}^n \ell(\tilde\pi_s^{(i)}) }
                { \suml_{i=1}^n m(\tilde\pi_s^{(i)}) - (n-1) }
        \le \f{  w_0  \suml_{i=1}^n m(\tilde\pi_s^{(i)}) + sd^{1/2} k} 
	      { \suml_{i=1}^n m(\tilde\pi_s^{(i)}) - (n-1) }
 \nonumber\\
       &\le&  w_0  
              + \f{ w_0 (n-1) + sd^{1/2}k}
	          {\suml_{i=1}^n m(\tilde\pi_s^{(i)}) - (n-1)} 
  \label{eq.wpns}		  
\end{eqnarray}
where 
$k = \#\{ i: w(\pi_s^{(i)}) >  w_0  $ .
By definition of $\tilde\pi_s$ either $m(\tilde\pi_s)\ge sd^{1/2}/w_0$,
or $m(\tilde\pi_s) = 1$, so that
\begin{eqnarray}
   m(\pi_{sn}) = 
   \suml_{i=1}^n m(\tilde\pi_s^{(i)}) - (n-1) 
&\ge& (n-k)( sd^{1/2}/w_0) + k - (n-1)
\nonumber \\ &
   \ge & (n-k)( sd^{1/2}/w_0 - 1 ).
\label{eq.msx}
\end{eqnarray}
Since $k/n\to P(w(\pi_s)>w_0)$ a.s. as $n\to\infty$,
combining \eqref{wpns} with \eqref{msx} gives \eqref{concat}
and the second inequality of the lemma follows from \eqref{msx}.
\eop

\begin{Proposition}
\label{Lbc}
$\beta = c(0+)$.
\end{Proposition}
\proof
Given a cycle through $m$ points in $[0,s]^d$, one can
make a path from $(0,\ldots,0)$ to $(s,\ldots,s)$
through these $m$ points using extra length at most
$2 s d^{1/2}$.
So for any $\delta >0$ we have, setting $s = n^{1/d}$
and $m = \lceil \delta n \rceil$,
\[ W_{n^{1/d}} \leq
\frac{L_n(\delta) + 2n^{1/d}d^{1/2}}{\lceil \delta n \rceil + 1}
 . \]
Letting $n \to \infty$ 
and using 
Proposition \ref{C1} 
and Theorem \ref{T1}(a),
we
see $\beta \leq c(\delta)$.
So $\beta \leq c(0+)$.

For the converse, consider the cube $[0,Ks]^d$ for even $K$,
partitioned into $K^d$ equal subcubes. 
In each subcube choose a diagonal path $\pi_s^{(i)}$
attaining the minimum $W_s$, 
so that concatenating these paths forms a cycle through the
large cube (see Figure 1).  

\putfigure{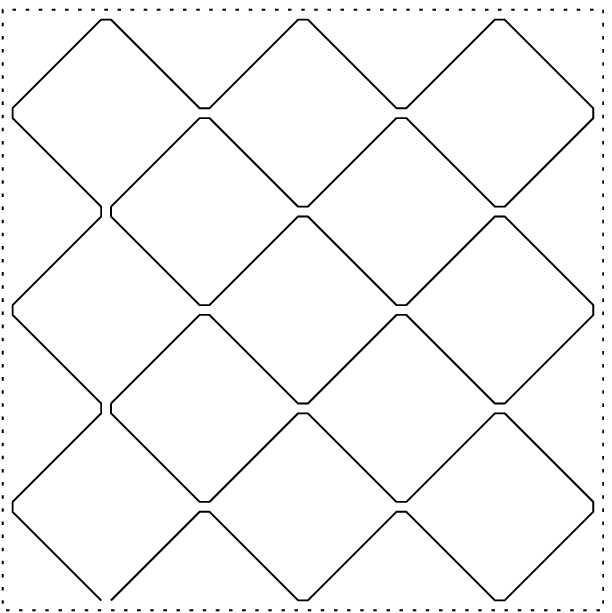}{A cycle built from subcube diagonals in the 6x6 cube}

Fix $\eps > 0$.
By Proposition \ref{C1}, 
for $s=s(\eps)$ large enough  
the random paths $\pi^{(i)}_s$  
satisfy
\[ P( w(\pi_s^{(i)}) > \beta + \eps ) \leq \eps . \]
Then by Lemma \ref{Lconcat}
applied with $w_0 = \beta + \eps$
we can construct cycles 
$\pi_{K^d}$ in $[0,Ks]^d$ such that
\begin{equation}
 P \left(
w(\pi_{K^d}) > \beta + \eps
+ \frac{\beta + \eps + sd^{1/2}\eps}{(1-\eps)(sd^{1/2}/(\beta + \eps)-1)}
+ \eps \right)
\to 0 \mbox{ as } K \to \infty .\label{wpb}
\end{equation}
By taking $s$ sufficiently large this simplifies to
\begin{equation}
 P(
w(\pi_{K^d}) > \beta + 2 \eps + 2 \beta \eps
) \to 0 \mbox{ as } K \to \infty .\label{wpb2}
\end{equation}
Also from Lemma \ref{concat}
there exists $\delta_s > 0$ such that
\[ P( \sfrac{m(\pi_{K^d})}{s^dK^d} \geq \delta_s)
\to 1 
\mbox{ as } K \to \infty 
. \]
On this event we have, for $\delta < \delta_s$,
\[ L(Ks, \delta) \le w(\pi_{K^d}) m(\pi_{K^d}). \]
So as $K \to \infty$
\[ P(L(Ks,\delta) \leq (\beta + 2 \eps + 2 \beta  \eps)s^dK^d \delta_s)
\to 1 \]
and then by
Lemma \ref{L2}
\[ \delta c(\delta) \leq (\beta + 2 \eps + 2 \beta \eps) \delta_s. \]
Letting $\delta\uparrow\delta_s$ then implies
\[ c(0+) \leq \beta+2 \eps + 2 \beta \eps, \]
so that letting $\eps \downarrow 0$ gives the desired inequality
\[ c(0+) \leq \beta. \]
%
%

Propositions \ref{C1} and \ref{Lbc} establish
part (c) of Theorem \ref{T1}.
For later use we record a small variation.


\begin{Lemma}[linear diagonal percolation]\label{ldp1}
For $\eta > 0$ let 
\[ W_s^{(\eta)} := \min_{\pi \in \PI_{0,s}: m(\pi) \leq \eta s}
   \frac{\ell(\pi)}{m(\pi)} \]
be the minimum average edge-length over paths 
across the diagonal of $[0,s]^d$ through
at most $\eta s$ points.
Then there exists a function $\beta(\eta)$ such that
\[ \lim_{s\to\infty} W_s^{(\eta)} = \beta(\eta) \quad\mbox{a.s.}, \]
\[ \beta(\eta) \downarrow c(0+) \mbox{ as } \eta \uparrow \infty . \]
\end{Lemma}

\proof 
The subadditivity argument used in Proposition \ref{C1} applies unchanged to 
$W_s^{(\eta)}$, giving the first limit.
The function
$\eta \to \beta(\eta)$ is {\em a priori} non-increasing,
and
\[ \lim_{\eta \to \infty} \beta(\eta) \geq \beta = c(0+) \]
by Proposition \ref{Lbc}.

Given $s < \infty$
let
\[ \mbox{ $\pi_s$ be a path attaining $W_s$
}\]
\[ \mbox{ $\pi^\prime_s$ be a path attaining $W^{(\eta)}_s$
for $\eta = 2s^{d-1}$ } . \]
Given $\eps > 0$ we can choose $s = s(\eps) < \infty$
sufficiently large that
\[ P(\pi^\prime_s \neq \pi_s) < \eps \]
(because $\eta s = 2s^d$ will likely exceed the 
Poisson$(s^d)$ number of points in $[0,s]^d$)
and
\[ P(w(\pi_s) > c(0+) + \eps) < \eps \]
by Propositions \ref{C1} and \ref{Lbc}.
So
\[ P(w(\pi^\prime_s) > c(0+) + \eps) < 2\eps .\]
Applying Lemma \ref{concat} to $\pi^\prime_s$ and
$w_0 = c(0+) + \eps$
gives the upper bound in
\[ \beta(2s^{d-1})
\leq \limsup_k w(\pi^\prime_{sk})
\leq c(0+) + \eps
+ \frac{c(0+) + \eps + sd^{1/2}2\eps}{
(\frac{sd^{1/2}}{c(0+)+\eps} - 1)
(1-2\eps)}  \]
and the lower bound follows from the construction in Lemma \ref{concat}.
Taking $\eps \to 0$, so that $s = s(\eps) \to \infty$,
\[ \lim_{\eta \to 0} \beta(\eta) \leq c(0+). \]



\subsection{Translation invariant distributions in the infinite model} \label{sec-tid}
Here we prove part (d) of Theorem \ref{T1} by relating
Approach 4 to Approach 3.
Fix $0 < \delta \leq 1$ and recall
$(\xi_i)$ denotes a Poisson point process of rate $1$ per unit volume in $\Reals^d$.
Recall from section \ref{sec-sparse}
the definition of $L(s,\delta)$.
Write $\EE(s,\delta)$ and $\VV(s,\delta)$ for the edge-set
and the vertex-set of a cycle attaining length 
$L(s,\delta)$.
For a cube $\Cube$, vertex-set $\VV$ and edge-set $\EE$,
write
$|\VV \cap \Cube|$ for the number of vertices of $\VV$ in $\Cube$,
and write
$\len(\EE \cap \Cube)$
for the total length of edges of $\EE$ restricted to $\Cube$.

Write $U_s$ for a uniform random position in $[0,s]^d$,
and write
$ \tilde{\xi}^{(s)}_i = \xi_i - U_s $
for the positions of the Poisson points relative to the
``random origin" $U_s$;
then write
$\widetilde{\EE}(s,\delta)$
for the corresponding set of relative positions of edges
of the tour attaining $L(s,\delta)$:
\[ 
\widetilde{\EE}(s,\delta)
= \{ (\tilde{\xi}^{(s)}_i,\tilde{\xi}^{(s)}_j): \ 
(\xi_i,\xi_j) \in \EE(s,\delta) \} . \]
Now the pair
$\left( 
(\tilde{\xi}^{(s)}_i), 
\widetilde{\EE}(s,\delta)
\right)$
takes values in the space
$\mathbf{S}^*$ 
of point-sets and paths defined as the space
$\mathbf{S}$
in section \ref{sec-Equiv},
except that for
$\mathbf{S}^*$ 
we allow cycles in addition to doubly-infinite paths.

There is a natural metric topology on
$\mathbf{S}^*$ 
obtained by regarding it as a space of marked point processes.
For each $s$ the point process 
$(\tilde{\xi}^{(s)}_i)$ 
is exactly a Poisson process.
By letting $s \to \infty$ through some subsequence
we can define a limit
\[
\left( 
(\tilde{\xi}^{(s)}_i), 
\widetilde{\EE}(s,\delta)
\right)
\cd \left( (\xi_i), \EE(\delta) \right) 
\mbox{ on } 
\mathbf{S}^* 
\]
where $\EE(\delta)$ is an edge set on some subset $\VV(\delta)$ of vertices
of the Poisson point process $(\xi_i)$.
It is easy to check that the edge-set forms doubly-infinite paths
(rather than finite cycles)
and so the right side has some distribution $\mu$ on 
$\mathbf{S}$.
From the uniform distribution for $U_s$ 
it is easy to check that
$\mu$ is translation invariant.
We will show (in the notation of section \ref{sec-Equiv})
\[ 
\delta(\mu) = \delta; \quad 
\ell(\mu) \leq c(\delta)
\]
which immediately implies
$\bar{c}(\delta) \leq c(\delta)$.
With a little extra effort we could prove 
$\ell(\mu) = c(\delta)$, but this doesn't help.

From the fact that $\VV(s,\delta)$
contains exactly $\lceil \delta N(s) \rceil$
of the $N(s)$ Poisson points in $[0,s]^d$,
it is clear that 
for fixed $r > 0$
\[
\lim_{s \to \infty}
E \left| 
\VV(s,\delta) \cap (U_s + \Cube_r)
\right| 
= \delta r^d \]
and therefore
\[
E \left| 
\VV(\delta) \cap \Cube_r 
\right| 
= \delta r^d ,\]
the interchange of limits being justified by the fact that
the total number of Poisson points in 
$U_s + \Cube_r$ has fixed distribution.
This tells us that $\delta(\mu) = \delta$.
Now fix $b < \infty$ and let 
$\EE^b(\delta)$ and $\EE^b(s,\delta)$
be the subsets of
$\EE(\delta)$ and $\EE(s,\delta)$
obtained by taking only edges of length $\leq b$.
By translation invariance
\[ E \len(\EE^b(\delta) \cap \Cube_r)
= \delta G(b) r^d \]
for some $G(b)$, and
then
\[ \ell(\mu) = \delta  \int_0^\infty b dG(b) . \]
By weak convergence,
letting $s \to \infty$ through a subsequence,
\begin{eqnarray*}
 E \len(\EE^b(\delta) \cap \Cube_1)
&=& 
\lim_s E \len(\EE^b(s,\delta) \cap (U_s + \Cube_1))\\
&=&
\lim_s s^{-d} E \len(\EE^b(s,\delta)).
\end{eqnarray*}
Since
$\delta c(\delta) = 
\lim_s s^{-d} E \len(\EE(s,\delta))$,
we can let $b \to \infty$ and use Fatou's lemma
to conclude
$\ell(\mu) \leq c(\delta)$.

For the converse, let $\mu$ be a translation invariant
distribution attaining $\bar{c}(\delta)$.
That is, $\mu$ specifies an edge-set $\EE$ on a vertex-subset
$\VV \subset (\xi_i)$ such that\\
(i) 
$E|\VV \cap \Cube_r| = \delta r^d \quad \forall r > 0$;\\
(ii)
$E \len(\EE \cap \Cube_r) = \delta \bar{c}(\delta) r^d \quad \forall r > 0$.\\
Now consider large $s$ and small $\eta > 0$.
The intersection of $\EE$ and $\Cube_s$
consists of a set of paths, each of which enters
$\Cube_s$ at some point on some face, and exits at some
point on some face.
Consider the subset of paths which intersect
$[\eta s, (1-\eta)s]^d$, and
write 
$\EE^{s,\eta}$
for the edges in this subset,
truncating an edge which crosses a face at the crossing point.
We extend this edge-set into a cycle in $\Cube_s$ as follows.
Suppose
the number of crossing points on each face is some even number 
(otherwise add superfluous edges, making no asymptotic difference);
let $D$ be the total number of crossing points.
On each face create a tour (in the face) of the crossing points;
by Lemma \ref{L1}(a) this has length at most
$A_1 s D^{(d-2)/(d-1)}$.
Within each such tour replace alternate edges by double edges.
The collection of these within-face edges, and the paths through
the interior of $\Cube_s$, form a connected graph where each
vertex has even degree, so we can find a Eulerian cycle; 
write 
$\FF^{s,\eta}$
for the edges in this cycle.
So
\begin{equation}
 E \len(\FF^{s,\eta}) 
\leq \delta \bar{c}(\delta) s^d
+ (4d)A_1s \left(ED \right)^{(d-2)/(d-1)} . \label{FAE}
\end{equation}
Because each crossing point is associated with some path-segment
of length $\geq \eta s$ from the face of $\Cube_s$
to the face of 
$[\eta s, (1-\eta)s]^d$, 
we have
\[ D \eta s \leq \len(\EE \cap \Cube_s) \]
and so by (ii)
\[ ED \leq (\eta s)^{-1} \delta \bar{c}(\delta) s^d . \]
Substituting into (\ref{FAE}) gives
\begin{equation}
 \limsup_s s^{-d} 
 E \len(\FF^{s,\eta}) 
\leq \delta \bar{c}(\delta) . \label{sEL}
\end{equation}
Now the cycle
$\FF^{s,\eta}$
goes through all the vertices
$\VV \cap 
[\eta s, (1-\eta)s]^d$.
By (i)
\[ E \left|
\VV \cap 
[\eta s, (1-\eta)s]^d
\right| = 
\delta (1-2\eta)^ds^d
\]
and by comparison with the total number of Poisson points in $\Cube_s$
we get a crude bound
\[ \var \left|
\VV \cap 
[\eta s, (1-\eta)s]^d
\right| \leq 
s^d(s^d + 1)
\leq 2s^{2d} \quad (s \geq 1) . \]
Now consider $k \geq 2$.
Using geometric subadditivity
(Lemma \ref{L1}(c))
and independent copies of $\FF^{s,\eta}$ on each of the $k^d$ subcubes,
we can construct a cycle
$\GG^{ks,\eta}$
in $[0,ks]^d$
such that (\ref{sEL}) extends to
\begin{equation}
 \limsup_s (ks)^{-d} 
 E \len(\GG^{ks,\eta}) 
\leq \delta \bar{c}(\delta) . \label{sEL2}
\end{equation}
Chebyshev's inequality implies
\begin{eqnarray*}
 P (
\mbox{ number of vertices in }
\GG^{ks,\eta} &\leq& 
\delta (ks)^d ((1-2\eta)^d - \eta)
\\&
\leq& \frac{2k^ds^{2d}}{(\delta k^ds^d\eta)^2}
= 2\delta^{-2}\eta^{-2} k^{-d} 
. \end{eqnarray*}
Taking $k = k(s) \to \infty$ sufficiently slowly, 
(\ref{sEL2}) remains true, so by definition of $c(\cdot)$
and the uniform boundedness property
\[
\delta ((1-2\eta)^d - 2 \eta)
\ c(
\delta ((1-2\eta)^d - 2 \eta)
) \leq \delta \bar{c}(\delta) . \]
Letting $\eta \to 0$ and using continuity of $c(\cdot)$
shows 
$c(\delta) \leq \bar{c}(\delta)$.

\subsection{The limit for $T_m$} \label{sec-Tm}
Recall the definition (Approach 1) of $T_m$.
If we could prove a sublinear growth property, that
the optimal path stays within a ball of radius $o(m)$,
then part (e) of Theorem \ref{T1} would follow easily from the other parts.
Because we cannot prove this property we use a more circuitous
method which eventually 
(Lemma \ref{Lfp})
compares a general path from the origin with either
a cycle through the origin
(Lemma \ref{fl1})
or a path between diagonals of a cube
(Lemma \ref{ldp3}).

Recall Lemma \ref{ldp1} on linear diagonal percolation.
We need a stronger result.
Given a realization
of the Poisson point process in $[0,s]^d$, first remove some $\mu s$ 
points, then look for the optimal diagonal path through at most $\eta s$ 
of the remaining points.  
We show that the same constant $\beta(\eta)$ 
appears in the limit.
\begin{Lemma}\label{ldp2}
Given the points of the Poisson process in $[0,s]^d$, 
for any $\eta > 0$, $\mu>0$ let 
\[ W_s^{(\eta,\mu)} := \max_{\Delta} \min_{\pi: \pi\cap\Delta=\emptyset}
			 \f{ \ell(\pi) }{ m(\pi) }
\]
where the maximum is taken over all subsets $\Delta$ of points with size at most $\mu s$
and the minimum over all paths $\pi \in \PI_{0,s}$ through at most $\eta s$ points that exclude $\Delta$.
Then 
\[ \lim_{s\to\infty} W_s^{(\eta,\mu)} = \beta(\eta) \mbox{ a.s.} \]
\end{Lemma}

\proof We consider the two-dimensional case.
The same proof works for dimension $d\ge 3$ without significant changes.

Consider the square
$[0,k^2s]^2$, 
partitioned into subsquares of side $s$.
Consider the $2k+1$ subsquares
$(S_i, \ 1 \leq i \leq 2k+1)$
whose lower left corners are
$(k,0), (k,1), (k-1,1), (k-1,2), \ldots, (1,k), (0,k)$.
For each $i$ consider the ``staircase" consisting of
$k^2 - k$ diagonally-adjacent subsquares starting with $S_i$.
Note that all these subsquares are distinct.
See Figure 2.


\putfigure{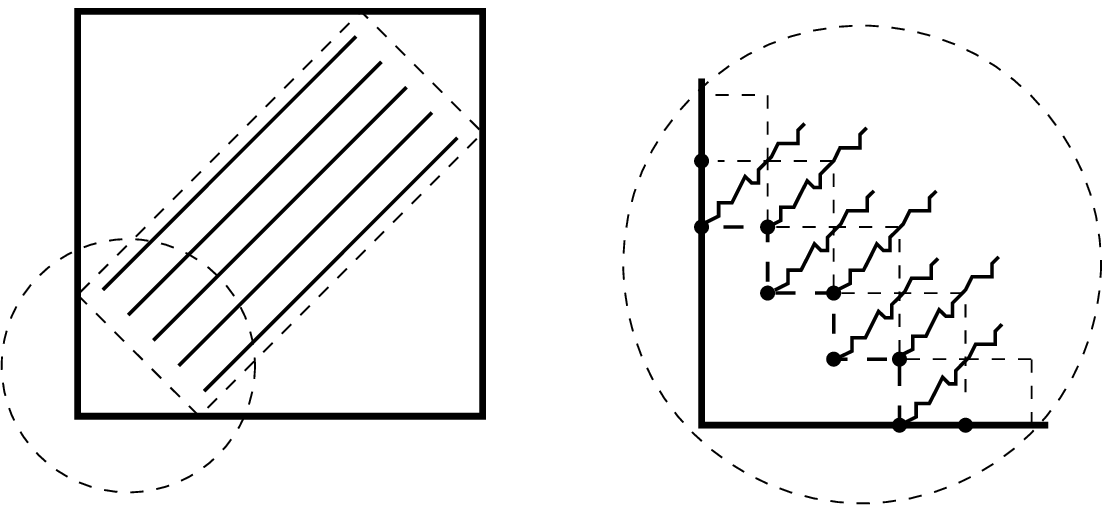}{$2k+1$ diagonal paths in a large square}

For each such staircase consider a joint path $\pi_k^{(i)}$,
obtained by concatenation of $k^2-k$ paths each attaining $W_s^{(\eta)}$
in its own subcube.

Given $\eta > 0$ and $\eps > 0$,
using Lemma~\ref{Lconcat} and Lemma~\ref{ldp1} we can choose large $s$ 
so that
\[ \limsup_k w(\pi^{(i)}_k) \le \beta(\eta)+\eps, \]
thus with probability $\to 1$ as $k\to\infty$ we can choose at least $3 k/2$
staircases $i$ satisfying
\begin{equation}\label{peb1}
m(\pi_k^{(i)}) \leq \eta k^2s, \qquad w(\pi_k^{(i)}) \leq \beta(\eta)+2\eps.
\end{equation}
%
For any subset of points $\Delta$ of size at most $\mu k^2$ there are at least $k$ 
staircases containing no more than $\mu k$ points of $\Delta$ each.
Thus with probability $1$, for all sufficiently large $k$
and all such $\Delta$ we can find some $\pi_k^{(i)}$ satisfying 
(\ref{peb1}) and
\[
\mbox{
$\pi_k^{(i)}$ contains $\leq \mu k$ points of $\Delta$}.
\]
From such a path $\pi_k^{(i)}$ we create a path
$\tilde\pi_k \in \PI_{0,k^2s}$ by
first deleting the points in $\Delta$ and 
then adding two edges from the path-ends to the corners
to $[0,k^2s]^d$.  
These changes have asymptotically negligible effect on
average edge-length, so we have shown
\[ \lim_{k\to\infty} W_{k^2s}^{(\eta,\mu)} 
\leq \beta(\eta) +  2\eps 
\mbox{ a.s. } \]
and the result follows easily.

\eop

Next we will use Lemma \ref{ldp2} to show that when one weakens conditions on $W_s^{(\eta)}$ 
by allowing paths through points outside the cube $[0,s]^d$,
this doesn't affect the limit $\beta(\eta)$.

\begin{Lemma}[unbounded linear percolation]\label{ldp3}
Let $\tilde W_s^{(\eta)}$ be the minimal average step length over all
paths from $(0,...0)$ to $(s,...s)$ through some points of a rate $1$
Poisson point process in $\Reals^d$, with at most $\eta s$ steps.
Then
\[ \lim_s \tilde W_s^{(\eta)} = \beta(\eta) \quad\mbox{in probability.}\]
\end{Lemma}

\proof 
Fix $\eta$ and consider the possibility that, for some $\beta'<\beta(\eta)$
\begin{equation}
 \limsup_s P( \tilde W_s^{(\eta)} < \beta' )  > 0. \label{pws}
 \end{equation}
 It is enough to show this cannot happen.
 
Let $\pi_s$ be the optimal path for $\tilde W_s^{(\eta)}$.
There exists a constant $\rho$ such that $\pi_s$ is contained
inside  $[-\rho s, s+\rho s ]^d$ a.s. for large $s$.
In the cubes $[-\rho s,0]^d$ and $[s, s+\rho s]^d$ 
take two diagonal paths $\pi_s^1$, $\pi_s^2$ with no more than $\eta\rho$
steps each, that attain minimal step length while avoiding the 
points of $\pi_s$.  By Lemma \ref{ldp2} the average edge-length in 
$\pi_s^1$ and in $\pi_s^2$ is asymptotically $\beta(\eta)$.
Then the concatenated path $\{\pi_s^1 , \pi_s,  \pi_s^2\}$ satisfies
the constraints in Lemma \ref{ldp1} for a
diagonal path across and within the large cube $[-\rho s, (\rho+1) s]$ 
with at most $\eta (1 + 2 \rho)s$ points. 
On the event in (\ref{pws}) this path has average step length 
less than some $\beta^{\prime \prime} <\beta(\eta)$.
So (\ref{pws}) contradicts the conclusion of Lemma \ref{ldp1}.
\eop

\begin{Lemma}[free cycle]\label{fl1}
Let $L_m(0)$ be the minimum length over all cycles on $m$ points
of the Poisson process that pass through the origin. Then
\[
\lim_{m\to\infty} L_m(0)/m \geq c(0+) \quad\mbox{in probability.}
\]
\end{Lemma}
{\em Remarks.} 
Such a cycle contains two edges from the origin to points
of the Poisson process.
Our results ultimately imply the inequality is really an equality,
but we won't need to prove that now.

\proof
Assume, to get a contradiction, that for some $c < c(0+)$
\begin{equation}
\limsup_{m\to\infty} P( L_m(0)/m < c ) > 0.
\label{Lm0}
\end{equation}
Let $\rho>0$ be a constant such that the cycle $\pi_m$ attaining length
$L_m(0)$ is contained
inside a box $[-\rho m, \rho m]^d$ a.s. for large $m$.

Fix $\eta>0$. 
By Lemma \ref{ldp2} it's possible to choose two paths $\pi_m^1$ and $\pi_m^2$ across
the diagonals of $[-\rho m,0]^d$ and $[0,\rho m]^d$ respectively in such a way that 
$\pi_m^1$, $\pi_m^2$ have no common points with $\pi_m$ except the origin, and
\[ n_1 \le \eta \rho m, \qquad l_1 = \beta(\eta) n_1 + o(m) \]
\[ n_2 \le \eta \rho m, \qquad l_2 = \beta(\eta) n_2 + o(m) \]
where $n_1$, $n_2$ and $l_1$, $l_2$ denote the number of edges and the
length of $\pi_m^1$ and $\pi_m^2$.

One can adjust a few edges around the origin to get a concatenated
path $\{\pi_m^1, \pi_m , \pi_m^2\}$ with length
$L^\prime \le l_1 + l_2 + L_m(0) $ and with $n^\prime = n_1+n_2+ m-2$ steps,
across the diagonal of the whole cube $[-\rho m, \rho m]$.
Observe
\begin{equation}
\f{ n_1 + n_2 + m }{2\rho m} \le \eta + \f1{2\rho }. \label{nnm}
\end{equation}
On the event 
$\{L_m(0)/m < c\}$ in (\ref{Lm0}) we have
\begin{eqnarray*}
 \f{ l_1 + l_2 + L_m(0) }{ n_1 + n_2 + m } - \beta(\eta) &\leq &
 \frac{o(m) + L_m(0) - \beta(\eta)m}{n_1 + n_2 + m}\\
 &\leq&
 \frac{o(m) + (c - \beta(\eta))m}{n_1 + n_2 + m}\\
 &\leq& \frac{c - \beta(\eta)}{1 + 2 \eta \rho}
+ o(1)
 \leq \frac{c - c(0+)}{1 + 2 \eta \rho}
+ o(1)
 \end{eqnarray*}
 the final two inequalities because
 $c - \beta(\eta) \leq c - c(0+) < 0$.
 Because of (\ref{nnm}) we can apply Lemma \ref{ldp1} and deduce
 
\[ \beta(\eta + \sfrac{1}{2\rho }) \le \beta(\eta) - \f{c(0+)-c}{2\eta \rho +1}. \]
But this implies $\lim_{\eta \to \infty} \beta(\eta) = -\infty$, which is impossible.
\eop

We can finally prove part (e) of Theorem \ref{T1}.
\begin{Lemma}[free path]
\label{Lfp}
Let $T_m$ be the minimal length over all paths through $m$ points 
of the Poisson process starting at the origin.
Then
\[
\lim_m T_m / m = c(0+) \quad\mbox{in probability.}
\]
\end{Lemma}
\proof
The upper bound
\[
\lim_k P( k^{-1}T_k \leq c(0+) + \eps) = 1 \ \forall \eps > 0
\]
was given at (\ref{4bound}).
For the lower bound suppose, for some $c<c(0+)$
\begin{equation}
\limsup_m P( T_m/m < c )  > 0.
\label{Tmm}
\end{equation}
Let $\pi_m$ be a path attaining length $T_m$ 
and let $x_m$ be the position of its end-point.
Consider
\[ r_* = \sup \{ r: \limsup_m P( |x_m| > r m, T_m < c m ) > 0 \}. \]
Clearly $r_* < \infty$.
Suppose first $r_* = 0$. 
Then we can construct a free cycle from $\pi_m$ by adding an edge from
$x_m$ to the origin, and such cycles have average length
$\leq c + o(1)$, contradicting Lemma \ref{fl1}.

Alternatively suppose $r_*>0$.  Then using rotational invariance 
of the Poisson process,
for arbitrary $\eps>0$ 
\[ \limsup_m P( x_m \in [(r-\eps)m,(r+\eps)m]^d \mbox{ and } T_m < c m ) >0. \]
Adding an edge from $x_m$ to $((r+\eps)m,\ldots,(r+\eps)m)$
gives a path between the diagonal corners of $[0,(r+\eps)m]^d$ 
through $m$ points
with average edge length at most $c+2\eps d^{1/2}$, which contradicts 
Lemma \ref{ldp3} when $\eps$ is small (because $\beta(\eta) \geq c(0+)$).  
So (\ref{Tmm}) is false and
we have proved
$T_m/m \to c(0+)$ in probability.
%
\eop

%
\section{Final Remarks}
\subsection{Subadditivity and cost-reward problems}
\label{sec-newsub}
The technique in section \ref{sec-Ws} seems applicable in many
contexts where subadditivity is used.
For instance, to modify the context of first-passage percolation on the lattice,
suppose that for each edge $e$ there is a ``cost" $c(e)$
and a ``reward" $r(e)$, so that for each path
$\pi = (e_1.e_2,\ldots,e_n)$ there is a cost and a reward
$(c(\pi),r(\pi)) = 
(\sum_i c(e_i),\sum_i r(e_i))$.
Then one can study
by the same technique
\[ A_n:= \min_\pi c(\pi)/r(\pi) \]
minimized over paths across the diagonal of $[0,n]^d$.
Developing a general theorem which establishes limits
$n^{-1}A_n \to a_0$ a.s.
in such settings 
would be a natural research project.

%


\begin{thebibliography}{10}

\bibitem{me109}
D.J. Aldous.
\newblock Percolation-like scaling exponents for minimal paths and trees in the
  stochastic mean-field model.
\newblock {\em Proc. R. Soc. Lond. Ser. A Math. Phys. Eng. Sci}, 461:825--838,
  2005.

\bibitem{BHH59}
J.~Beardwood, H.J. Halton, and J.M. Hammersley.
\newblock The shortest path through many points.
\newblock {\em Proc. Cambridge Phil. Soc.}, 55:299--327, 1959.

\bibitem{BKS03}
I.~Benjamini, G.~Kalai, and O.~Schramm.
\newblock First passage percolation has sublinear distance variance.
\newblock {\em Ann. Probab.}, 31:1970--1978, 2003.

\bibitem{durrett-kesten}
R.~Durrett.
\newblock {H}arry {K}esten's publications: a personal perspective.
\newblock In M.~Bramson and R.~Durrett, editors, {\em Perplexing Problems in
  Probability}, pages 1--33. Birkhauser, 1999.

\bibitem{gri99}
G.R. Grimmett.
\newblock {\em Percolation}.
\newblock Springer-Verlag, Berlin, 2nd edition, 1999.

\bibitem{howard-fpp}
C.~D. Howard.
\newblock Models of first-passage percolation.
\newblock In H.~Kesten, editor, {\em Probability on Discrete Structures},
  volume 110 of {\em Encyclopaedia of Mathematical Sciences}, pages 125--173.
  Springer-Verlag, 2004.

\bibitem{JRS04}
J.L. Jacobsen, N.~Read, and H.~Saleur.
\newblock The traveling salesman problem, conformal invariance, and dense
  polymers.
\newblock {\em Phys. Rev. Lett.}, 93:038701--038704, 2004.
\newblock arXiv:cond-mat/0403277.

\bibitem{kesten-FPP}
H.~Kesten.
\newblock First-passage percolation.
\newblock In {\em From Classical to Modern Probability}, number~54 in Progr.
  Probab., pages 93--143. Birkhauser, 2003.

\bibitem{meester-CP}
R.~Meester and R.~Roy.
\newblock {\em Continuum Percolation}.
\newblock Cambridge University Press, 1996.
\newblock Cambridge Tracts in Math. 119.

\bibitem{MP86}
M.~M\'{e}zard and G.~Parisi.
\newblock A replica analysis of the travelling salesman problem.
\newblock {\em J. Physique}, 47:1285--1296, 1986.

\bibitem{rhee-91}
W.~Rhee.
\newblock On the fluctuations of the traveling salesperson problem.
\newblock {\em Math. Oper. Res.}, 16:482--489, 1991.

\bibitem{ste81}
J.M. Steele.
\newblock Complete convergence of short paths and {K}arp's algorithm for the
  {T}{S}{P}.
\newblock {\em Math. Oper. Res.}, 6:374--378, 1981.

\bibitem{steele97}
J.M. Steele.
\newblock {\em Probability Theory and Combinatorial Optimization}.
\newblock Number~69 in CBMS-NSF Regional Conference Series in Applied Math.
  SIAM, 1997.

\bibitem{yukich-book}
J.E. Yukich.
\newblock {\em Probability Theory of Classical {E}uclidean Optimization
  Problems}.
\newblock Number 1675 in Lecture Notes in Math. Springer, 1998.

\end{thebibliography}

\def\cprime{$'$} \def\polhk#1{\setbox0=\hbox{#1}{\ooalign{\hidewidth
  \lower1.5ex\hbox{`}\hidewidth\crcr\unhbox0}}} \def\cprime{$'$}
  \def\cprime{$'$} \def\cprime{$'$}
  \def\polhk#1{\setbox0=\hbox{#1}{\ooalign{\hidewidth
  \lower1.5ex\hbox{`}\hidewidth\crcr\unhbox0}}} \def\cprime{$'$}
  \def\cprime{$'$} \def\polhk#1{\setbox0=\hbox{#1}{\ooalign{\hidewidth
  \lower1.5ex\hbox{`}\hidewidth\crcr\unhbox0}}} \def\cprime{$'$}
  \def\cprime{$'$} \def\cydot{\leavevmode\raise.4ex\hbox{.}} \def\cprime{$'$}
  \def\cprime{$'$} \def\cprime{$'$} \def\cprime{$'$} \def\cprime{$'$}
  \def\cprime{$'$} \def\cprime{$'$} \def\cprime{$'$} \def\cprime{$'$}
  \def\cprime{$'$}

\end{document}